\documentclass{ifacconf}
\usepackage{amsmath,amssymb}
\usepackage{mathrsfs,graphicx} 

\begin{document}
\begin{frontmatter}
\title{Recursive set-membership state estimation for linear non-causal time-variant differential-algebraic equation with continuous time}
\author[First]{Sergiy M. Zhuk} 

\address[First]{Taras Shevchenko Kyiv National University, Faculty of
  Cybernetics, Kyiv, Ukraine (e-mail: Serhiy.Zhuk@gmail.com).}
\begin{abstract} 
This paper describes a state estimation approach for non-causal time-varying linear descriptor equations with uncertain parameters. The uncertainty in the state equation and in the measurements is supposed to admit a 
set-membership description. The approach is based on the notion
of the \emph{linear minimax estimation}. Suboptimal minimax state estimation algorithm is introduced for DAEs with non-stationary rectangular matrices. Optimal algorithm is presented for DAEs with special structure of the matrices. A comparison of suboptimal and optimal algorithms is presented for 2D time-varying DAE with a singular matrix pencil. 
\end{abstract}
\end{frontmatter}
\section{Introduction}
Dynamical systems described by coupled differential and algebraic equations
(DAEs) arise naturally in many applications. In particular, DAEs occur in 
econometrics (singular dynamic Leotief systems~\cite{Luenberger1977}), 
modelling of constrained multibody systems~\cite{Mills2006}, 
electrical circuit synthesis~\cite{Reis2008}, 
bioprocess and chemical engineering~\cite{Becerra2001},
representing a common modelling tool~\cite{Mehrmann2005}. On can divide contributions to the theory \footnote{Solvability and numerical methods, control and observation for DAE and so on.} of DAEs into results for casual DAEs and results for non-causal systems. 
{\bf 1.1 Causal DAEs. } 
The solvability theory for finite-dimentional systems with constant coefficients  
\begin{equation}
  \label{eq:ldae}
  F x_t=Cx+Bf
\end{equation}
is based on the reduction of the matrix pencil $sF-C$ to the Kronecker canonical form: if $\mathrm{det}[sF-C]\not\equiv0$ then for all initial values $x(t_0)=x_0$ there exists a unique solution $x(\cdot)$. Changing the basis in the state space and differentiating exactly $s$
times ($d$ is an index of the pencil $sF-C$), one can reduce
\eqref{eq:ldae} to some equivalent Ordinary Differential Equation (ODE), provided $f(\cdot)$ is sufficiently smooth and meets some algebraic constraints. The details of the reduction process are presented in~\cite{Gantmacher1960}. The index $d$ of the pencil $sF-C$ is said~\cite{Campbell1983} to be an index of linear DAE \eqref{eq:ldae}. 

The notion of a Standard Canonical Form (SCF) allows to generalize the index approach on variable coefficients: for instance, in \cite{Campbell1983} it was shown that
\eqref{eq:ldae} with analitical $F(\cdot)$, $C(\cdot)$ and $B(\cdot)$ 
is solvable\footnote{DAE \eqref{eq:ldae} is said to be solvable at $[t_0,T]$ if
for every sufficiently smooth $f(\cdot)$ there exists at least one continuously
differentiable solution, provided $F(\cdot)$, $C(\cdot)$ are sufficiently smooth.} iff \eqref{eq:ldae} can be converted into SCF. 
In \cite{Campbell1987} it was noted
that not all solvable DAEs can be put into SCF and the solvable DAE is equal
to some differential-algebraic equation in the canonical form which generalize
SCF. In this regard, we say that \emph{the DAE is causal} if 1) it can be reduced (at least locally in the non-linear case) to normal ODE and 2) if it is solvable for the given intial condition $x_0$ and input $f(\cdot)$ then the solution is  unique. 

The geometry of
the reduction procedure for nonlinear causal DAEs $F(x,\dot x)=0$ was
investigeted in \cite{Reich1990,Rabier1994}, where the index of DAE was
defined as a smallest natural $d$ so that the sequence of the constraint
manifolds \cite{Reich1990} $$
M_k:=TW_{k-1}\cap M_{k-1}, 
W_k=\{x\in\mathbb R^n:(x,p)\in M_k\}
$$
with $M_0:=\{(x,p):F(x,p)=0\}$ becomes stationary for $k>d$. 
This coincides with the definition of the index of linear DAE. 
Solvability of the linear casual DAEs with impulses in the input was addressed in~\cite{Rabier1996}. Further discussion of the linear DAEs
solvability theory in finite-dimensions and related topics are presented in
\cite{Mehrmann2005,Samoilenko2000}.

Basic ideas of the index approach (reduction of the pencil $sF-C$ to the canonical form) were extended on systems with constant operator coefficients in \cite{Rutkas1975}, provided the poles of the operator-valued function $s\mapsto(sF-C)^{-1}$ are contained in some bounded vicinity of $0$. 

{\bf 1.2 Non-causal DAEs. } 
The non-causal DAE may have several solutions. For instance,
consider $$
\dfrac d{dt}\bigl[
\begin{smallmatrix}
  1&&0
\end{smallmatrix}\bigr](
\begin{smallmatrix}
  x_1\\x_2
\end{smallmatrix})=\bigl[
\begin{smallmatrix}
  c_1&&c_2
\end{smallmatrix}\bigr](
\begin{smallmatrix}
  x_1\\x_2
\end{smallmatrix})+f(t), \bigl[
\begin{smallmatrix}
  1&&0
\end{smallmatrix}\bigr](
\begin{smallmatrix}
  x_1\\x_2
\end{smallmatrix})(t_0)=x_0
$$  Let $x_2(\cdot)\in\mathbb
L_2(t_0,T)$,  $f(\cdot)\in\mathbb L_1(t_0,T)$ and
$x_0\in\mathbb R$. By definition put $$
x_1(t):=\exp(c_1(t-t_0))x_0+\int_{t_0}^t\exp(c_1(t-s))c_2x_2(s)+f(s)ds
$$ 
It is clear that any solution of this DAE is given by the formula $
t\mapsto(x_1(t),x_2(t))^T$. According to a behavioral approach
\cite{Ilchmann2005} one can think about $x_2$ as an input or as a part of the
system state. Reduction of the non-causal DAEs with variable matrices was studied in~\cite{Shlapak1975} and \cite{Eremenko1980}, provided $F(\cdot)$ has a constant rank. In this case DAE may be splitted into differential and algebraic equation. Further splitting is possible under some restrictions on matrices of resulting system. 



In this paper we focus on the inverse problem for linear non-causal DAEs with rectangular non-stationary matrices: given measurements $y(t)$, $t\in[t_0,T]$ of the DAE' solution $x(\cdot)$, to reconstruct $x(\cdot)$.  
Here $x(\cdot)$ is said to be a solution of   
\begin{equation}
  \label{eq:FDAE}
  (Fx)_t(t)=C(t)x(t)+f(t)
\end{equation} 
with initial condition $Fx(t_0)=0$, if $Fx(\cdot)$ is absolutely continuous
function, $(Fx)_t$ belong to $\mathbb L_2(t_0,T,\mathbb R^m)$, $x(\cdot)$ verifies \eqref{eq:FDAE} almost everywere and
$Fx(t_0)=0$ holds. This definition guarantees that a linear mapping $D$, induced by \eqref{eq:FDAE}, is closed as a mapping from $\mathbb L_2(t_0,T,\mathbb R^n)$ to $\mathbb L_2(t_0,T,\mathbb R^m)$ \cite{Zhuk2007}. This, in turn, allows to properly define the system~\cite{Zhuk2007}, adjoint to~\eqref{eq:FDAE}, that is of primary interest in optimal control. In addition, the method of matrix pensils is sometimes difficult to apply in the finite-dimensional optimal control theory. For instance \cite{Ozcaldiran1989}, applying the linear proportional feedback $f=Kx$ to \eqref{eq:ldae} with regular pencil $\mathrm det(sF-C)\ne0$ one may arrive to the system with singular pensil $\mathrm
det(sF-C-BK)\equiv 0$. Therefore it is reasonable to apply the above definition of the DAEs' solution in the control framework. A much general one is presented in \cite{Kurina2007}, where a properly
stated leading term $A(t)\dfrac d{dt}F(t)x$ is used in order to give a
feedback solution to LQ-control problem with DAE 
constraints. 

In what follows we assume that $y(t)=H(t)x(t)+\eta(t)$, $x(\cdot)$ is a solution of \eqref{eq:FDAE} in the sense of the above definition, the noise $\eta(\cdot)$ is a realization of a random process $\Psi$, the input $f(\cdot)$ is uncertain and belongs to the given set $G$. The aim is to construct a worst-case estimation of the inner product $\langle\ell,Fx(T)\rangle$, $\ell\in\mathbb R^n$ as a function of $y(\cdot)$, assuming that 
$G$ is bounded set and the correlation function of $\Psi$ belongs to the given
bounded set $R$ of matrix-valued functions. This problem was solved in~\cite{Nakonechnii1978}, provided $F=I$. 
The case of deterministic measurement's noise was addressed in \cite{Bertsekas1971,Tempo1985,Kurzhanski1997}, where the optimal worst-case estimation is shown to be a dynamical system, describing the evolution of the
 central point of the ODEs reachability set, consistent with observations. 

In this paper, we generalize the theory of minimax state estimation~\cite{Nakonechnii1978} on a class of linear non-causal DAEs: $F\in\mathbb R^{m\times n}$ and $t\mapsto C(t)\in\mathbb R^{m\times n}C(t)$ is continuous on $[t_0,T]$. The same results can be proved for the deterministic noise $\eta(\cdot)$ and ellipsoidal bound for uncertain $f(\cdot)$ and $\eta(\cdot)$, giving the generalization of \cite{Bertsekas1971,Kurzhanski1997}.    
The major contributions of this paper is an implementation of the abstract  Generalized Kalman Duality principle~\cite{Zhuk2009d} for non-causal time-dependent DAEs (Theorem~\ref{t:1}). Duality allows to find and exact expression for the worst-case estimation error and to establish the necessary and sufficient conditions on $\ell$ for the worst-case error to be finite. These conditions, in turn, defines some subspace $\mathcal L(T)$ in the state space of~\eqref{eq:FDAE}, which is called a minimax observable subspace. In fact, $\mathcal L(T)$ describes an "observable" (in the minimax sense) part of $x(T)$ with respect to the measured $y(t)$, $t\in[t_0,T]$: if $\ell\in\mathcal L(T)$ then we can provide the worst-case estimation of $\langle\ell,Fx(T)\rangle$ with finite worst-case error (which describes the measure of how poor the estimation quality may be); otherwise the state $x(T)$ is not observable in the direction $\ell$, that is for any estimation of $\langle\ell,Fx(T)\rangle$ the estimation error varies in $[0,+\infty]$, so that, for any linear estimation and natural $N$ there is a realization of uncertain parameters $f(\cdot)$ and $\eta(\cdot)$ such that the estimation error will be greater than $N$. Note, that the notion of the minimax observable subspace $\mathcal L(T)$ is a implementation (in the case of DAEs) of the abstract minimax observability concept, presented in \cite{Zhuk2009d}. Some aspects of classical observability for DAEs were considered in \cite{Campbell1991} for causal systems and in \cite{Frankowska1990} for non-causal systems. 

As a result of application of Generalized Kalman Duality for non-causal non-stationary DAE we derive a suboptimal worst-case state estimation algorithm (Corollary~\ref{c:eme}). The algorithm gives a suboptimal estimation of the projection of the state $x(T)$ onto the minimax observable subspace $\mathcal L(T)$. It is sequential, that is the algorithm is represented in terms of the unique solution to a Cauchy problem for some ODE, which has a realization of the observations $y(t)$, $t\in[t_0,T]$ as the input. Therefore, it is sufficient to know measurements $y(t)$, $t\in[T,T_1]$ and the estimation at $t=T$ in order to compute the estimation of $x(T_1)$. The algorithm works for "non-Gaussian noise" $\eta$ unlike the family of Kalman-like estimators. The optimal algorithm is also presented (Corollary~\ref{c:me}), provided the matrices of DAEs have "some regularity" (Proposition~\ref{p:2}). Kalman filtering approach was previously applied to linear DAEs with constant coefficients in \cite{Gerdin2007,Darouach1997}, provided $sF-C$ is regular. In this regard we note that the latter assumption can be substituted by the less restrictive one: $sF'-C'-H'H$ is regular (see Example above). Further information on Kalman filtering for causal DAEs is to be found at \cite{Xu2006}. A worst-case state estimation for non-causal linear continuous DAEs with non-stationary rectangular matrices was not considered in the literature before. The notion of the minimax observable subspace was applied in \cite{Zhuk2009c} in order to construct the optimal state estimation algorithm for discrete time non-causal DAEs.

\emph{Notation}: $\mathrm E\eta$ denotes the mean of the random element $\eta$, $\mathrm{int}\, G$ denotes the interior of $G$, 
$f(\cdot)$ or $f$ denotes some element of the functional
space, 
$f(t)$ denotes the
value of the function $f$ at time $t$, $\mathbb
L_2(t_0,T,\mathbb R^m)$ denotes the space of square-integrable functions with values in $\mathbb R^m$, $\mathbb H^1(t_0,T,\mathbb R^m)$ denotes the space of absolutely continuous functions with $L_2$-derivative and values in $\mathbb R^m$, the superscript $'$ denotes the operation of taking an adjoint, 
$c(G,\cdot)$ denotes the support function of some set $G$, $\delta(G,\cdot)$
denotes the indicator\footnote{$\delta(G,f)=0$ if $f\in G$ and $+\infty$
otherwise.} of $G$, 
$\mathrm{dom} f=\{x:f(x)<\infty\}$; 
$\langle\cdot,\cdot\rangle$ 
denotes the inner 
product in Hilbert space $\mathcal H$, $\mathbb R^n$ denotes $n$-dimensional
Euclidean space over real field, 
$S>0$ means $\langle Sx,x\rangle>0$ for all $x$ from within appropriate Hilbert space, $L'$ denotes adjoint operator, $F'$ denotes transposed matrix, $F^+$ denotes pseudoinverse matrix.  
\section{Linear minimax estimation for DAEs}\label{s:obs}       
Consider a pair of systems
  \begin{equation}
    \label{eq:dae}
    \begin{split}
      & (Fx)_t(t)=C(t)x(t)+f(t),Fx(t_0)=0,\\
      & y(t)=H(t)x(t)+\eta(t),t\in[t_0,T],
    \end{split}
   \end{equation}
where $x(t)\in\mathbb R^n$, $f(t)\in\mathbb R^m$, $y(t)\in\mathbb R^p$,
$\eta(t)\in\mathbb R^p$ represent the state, input, observation and
observation's noise respectively, $F\in\mathbb R^{m\times n}$, $f(\cdot)\in\mathbb L_2(t_0,T)$, $C(t)$ and $H(t)$ are continuous matrix-valued functions, $t_0,T\in\mathbb R$.  

According to~\cite{Zhuk2007} we say that $x(\cdot)$ is a
solution of~\eqref{eq:dae} if $Fx(\cdot)\in\mathbb H^1(t_0,T,\mathbb R^m)$ and the derivative $(Fx)_t$ of $Fx(\cdot)$ coincides with the right side
of~\eqref{eq:dae} almost everywhere (a.e.) on $[t_0,T]$ and $Fx(t_0)=0$. 
\begin{rem}
  As $F\in\mathbb R^{m\times n}$ the pencil $F-\lambda C(t)$ is irregular~\cite{Muller} implying~\eqref{eq:dae} may have non-unique solution. In what follows we will refer such DAEs as non-causal~\cite{Muller}. 
\end{rem}
In the sequel we assume $\eta(\cdot)$ is
a realization of a random process $\Psi$ such that $E\Psi=0$ and 
\begin{equation}
  \label{eq:eta_bounds}
\Psi\in W=\{\Psi:E\int_{t_0}^T(R(t)\eta(t),\eta(t))\le 1\}
\end{equation}
and
\begin{equation}
  \label{eq:G}
  f(\cdot)\in G:=\{f(\cdot):\int_{t_0}^T(Q(t)f(t),f(t))\le 1\},
\end{equation}
where $Q(t)\in\mathbb R^{m\times m}$, $Q=Q'>0$, $R(t)\in\mathbb
R^{p\times p}$, $R'=R>0$ and $Q(t),R(t)$, $R^{-1}(t)$, $Q^{-1}(t)$ are
continuous functions of $t$ on $[t_0,T]$. 

Suppose $y(t)$ is observed in \eqref{eq:dae} for some $x(\cdot)$, $f\in G$ and
$\Psi\in W$. The purpose of this section is to construct an algorithm with the following property: given a realization $y(t)$, $t\in [t_0,T]$ of the random process $Y=Hx(\cdot)+\Psi$, the algorithm produces  
an estimation of a linear function $$
x(\cdot)\mapsto \langle\ell, Fx(T)\rangle
$$ having minimum mean-squared worst-case estimation error. 
In what follows we will refer this algorithm as an a priori minimax mean-squared estimation in the direction $\ell$ ($\ell$-estimation). 
Taking into account linearity of~\eqref{eq:dae} we will look for $\ell$-estimation among linear functions of $y(\cdot)$. Let us summarize the above discussion by rigorous mathematical definitions.   
\begin{defn}
Given $u(\cdot)\in\mathbb L_2(t_0,T,\mathbb R^p)$ and $\ell\in\mathbb R^m$ define a mean-squared worst-case estimation error\footnote{Here the $\sup_{\eta(\cdot)}$ means that we take the upper bound over all random processes $\Psi$ such that $E\int_{t_0}^T(R(t)\eta(t),\eta(t))\le 1$ for all realizations $\eta(\cdot)$ of $\Psi$. } 
\begin{equation}
  \begin{split}
    &\sigma(T,u,\ell):=\sup_{x(\cdot),f(\cdot),\eta(\cdot)}\{
E[\langle\ell, Fx(T)\rangle-u(y)]^2:\\
&(Fx)_t=Cx+f,Fx(t_0)=0,f(\cdot)\in G,\eta(\cdot)\in W\}
  \end{split}
\end{equation}
A function $\hat u(y)=\int_{t_0}^{T}(\hat u(t),y(t))dt$ is called
an a priori minimax mean-squared estimation in the direction $\ell$ ($\ell$-estimation) if $\inf_u\sigma(T,u,\ell)=\sigma(T,\hat u,\ell)$. 
The number
$\hat\sigma(T,\ell)=\sigma(\hat u,\ell)$ is called a minimax mean-squared a priori error in the direction $\ell$ at time-instant $T$ ($\ell$-error). 
The set $\mathcal L(T)=\{\ell\in\mathbb R^n:\hat\sigma(\ell)<+\infty\}$ is called a minimax observable subspace. 
\end{defn}
\subsection{Generalized Kalman Duality Principle}\label{s:gkd}
  The definition of the $\ell$-estimation and error generalizes the notion of the linear minimax a priori mean-squared estimation, introduced in \cite{Nakonechnii1978}. In order to find the $\ell$-estimation we will follow a common way of deriving the minimax 
  estimation \cite{Nakonechnii1978}: 
first step is to obtain the expression for the worst-case error by means of the suitable duality concept, that is to formulate a dual control problem; 
next step is to solve it and to derive the minimax estimation.  

Next theorem generalizes the celebrated Kalman duality
principle~\cite{Brammer1989} to non-causal DAEs. 
\begin{thm}[Generalized Kalman duality]
\label{t:1}
The $\ell$-error
is finite iff  
\begin{equation}
    \label{eq:zul}
    (F'z)_t(t)=-C'(t)z(t)+H'(t)u(t), F'z(T)=F'\ell
  \end{equation}
\eqref{eq:zul} has a solution $z(\cdot)$. In this case 
the problem $\sigma(u)\to\inf_u$ is equal to the following optimal control problem 
\begin{equation}
  \label{eq:umin}
  \begin{split}
 I(u)=&\min_v\{\int_{t_0}^T(Q^{-1}(z-v),z-v)dt\}\\
&+\int_{t_0}^T(R^{-1}u,u)dt\to\min_u,   
  \end{split}
 \end{equation}
provided $z(\cdot)$ obeys \eqref{eq:zul} and $v(\cdot)$ obeys homogeneous \eqref{eq:zul}. 
\end{thm}
An obvious corollary of the Theorem~\ref{t:1} is an expression for the minimax observable subspace $$
\mathcal L(T)=\{\ell\in\mathbb R^n:\exists u(\cdot),z(\cdot)\}\,(F'z)_t+C'z-H'u=0,\,F'z(T)=F'\ell\}
$$ 
\begin{pf}
Take $\ell\in\mathbb R^n$ and $u(\cdot)\in\mathbb L_2(t_0,T,\mathbb R^p)$ and suppose $\ell$-error is finite.  
There exists some $w(\cdot)\in\mathbb L_2(t_0,T,\mathbb R^m)$ so that $F'w(\cdot)\in\mathbb H^1(t_0,T,\mathbb R^n)$ and $F'w(T)=F'\ell$. A trivial example is $w(t)\equiv\ell$. It was proved in~\cite{Zhuk2007} that 
\begin{equation}
  \label{eq:ibp}
  \begin{split}
&\langle F'w(T),F^+Fw(T)\rangle-\langle Fx(t_0),F^+Fw(t_0)\rangle =    \\
&\int_{t_0}^T \langle(Fx)_t,w\rangle+\langle(F'w)_t,x\rangle dt 
  \end{split}
\end{equation}
if $Fx(\cdot)\in\mathbb H^1(t_0,T,\mathbb R^m)$ and $F'w(\cdot)\in\mathbb
H^1(t_0,T,\mathbb R^n)$. 
Noting that $F=FF^+F$ and using~\eqref{eq:ibp} and~\eqref{eq:dae} one derives  
\begin{equation}\label{eq:i1}
  \begin{split}
    \langle\ell,Fx(T)\rangle&=\langle F'\ell,F^+Fx(T)\rangle = \langle F'w(T),F^+Fx(T)\rangle\\ 
&=\int_{t_0}^T \langle(Fx)_t,w\rangle+\langle(F'w)_t,x\rangle dt\\
&=\int_{t_0}^T \langle f,w\rangle+\langle(F'w)_t+C'w,x\rangle dt,
  \end{split}
\end{equation}
Combining~\eqref{eq:i1} with $E\eta=0$ we have 
\begin{equation}
  \label{eq:i2}
\begin{split}
  E[\langle\ell,& Fx(T)\rangle-u(y)]^2\\
&=
  [\langle\ell, Fx(T)\rangle-\int_{t_0}^T\langle H'u,x\rangle dt]^2+
E[\int_{t_0}^T\langle u(t),\eta(t)\rangle dt]^2\\
&=[\int_{t_0}^T \langle f,w\rangle+\langle(F'w)_t+C'w-H'u,x\rangle dt]^2\\
&+E[\int_{t_0}^T\langle u(t),\eta(t)\rangle dt]^2
  \end{split}
\end{equation}
Combining~\eqref{eq:eta_bounds} with Cauchy inequality we obtain 
\begin{equation}
  \label{eq:etaSup}
\sup_\eta E(\int_{t_0}^{T}(u,\eta)dt)^2=\int_{t_0}^{T}(R^{-1}u,u)dt  
\end{equation}
\eqref{eq:etaSup} and $\sigma(u)<+\infty$ imply the third line of~\eqref{eq:i2} is bounded. Noting that $\int_{t_0}^T\langle f,w\rangle dt$ is bounded independently of $x(\cdot)$ one derives  
\begin{equation}
\label{eq:i3}
\sup_{x(\cdot)}
\{\int_{t_0}^T\langle(F'w)_t+C'w-H'u,x\rangle dt
:(Fx)_t=Cx+f,
f(\cdot)\in G\}<+\infty
\end{equation} 
It was proved in~\cite{Zhuk2009d} that 
\begin{equation}
  \label{eq:LfiG}
  \sup_{x\in\mathscr D(L)}\{\langle\mathcal L, x\rangle, Lx\in G\}=
\inf_{b\in\mathscr D(L')}\{c(G,b),L'b=\mathcal L\}
\end{equation} provided $\mathscr D(L):=\{x(\cdot)\in\mathbb L_2(t_0,T,\mathbb R^n):Fx(\cdot)\in\mathbb
H^1(t_0,T,\mathbb R^n),Fx(t_0)=0\}$ and   
\begin{equation}
  \label{eq:Lx}
  (Lx)(t)=(Fx)_t(t)-C(t)x(t),x(\cdot)\in\mathscr D(L)
\end{equation} 
It was proved in~\cite{Zhuk2007} that $\mathscr D(L'):=\{b\in\mathbb L_2(t_0,T,
\mathbb R^m):F'b(\cdot)\in\mathbb H^1(t_0,T,\mathbb R^m),F'b(T)=0\}$ and 
\begin{equation}
  \label{eq:sLz}
  L'b(t)=-(F'b)_t-C'(t)b(t),
b(\cdot)\in\mathscr D(L')
\end{equation}
provided $L$ is defined by~\eqref{eq:Lx}. Setting $\mathcal L:=(F'w)_t+C'w-H'u$ we see from~\eqref{eq:i3} that the right-hand part of~\eqref{eq:LfiG} is finite. Using~\eqref{eq:sLz} one derives  
\begin{equation}
  \label{eq:4}
  \inf\{c(G,b),-(F'b)_t-C'(t)b(t)=(F'w)_t+C'w-H'u\}<+\infty
\end{equation}
Thus 
\eqref{eq:4} implies $$
(F'z)_t+C'z=H'u,F'z=F'\ell
$$ with $z:=(w+b)$, $b(\cdot)\in\mathscr D(L')$. This proves~\eqref{eq:zul} has a solution $z(\cdot)$. Using integration-by-parts formula~\eqref{eq:ibp} and $E\Psi=0$ and \eqref{eq:etaSup} one derives easily 
\begin{equation}
  \label{eq:gerr}
  \sigma(u)=\sup_{f\in G_1}\int_{t_0}^T\langle w,f\rangle dt
+\int_{t_0}^{T}(R^{-1}u,u)dt 
\end{equation}
with $G_1$ denoting all $f(\cdot)\in G$ such that~\eqref{eq:dae} has a solution $x(\cdot)$.  
 
On the contrary, if $z(\cdot)$ is some solution of~\eqref{eq:zul} then one derives~\eqref{eq:gerr} as it has been already done above. Therefore, there are only two cases: $\ell$-error is infinite or \eqref{eq:gerr} holds.    

Note that
\begin{equation}
  \label{eq:supG1}
  \sup_{f\in G_1}(\int_{t_0}^T(f,z)dt)^2=\sup\{\langle f,z \rangle, f\in G\cap R(L)\}^2
\end{equation}
where $R(L)$ is the range of the linear mapping $L$
defined above by the rule~\eqref{eq:Lx}.
It was proved in~\cite{Zhuk2009d} that 
\begin{equation}
  \label{eq:supGRL}
  \sup\{\langle f,z \rangle, f\in G\cap R(L)\}=
\inf\{c(G,z-v), v\in N(L')\}
\end{equation}
provided $\mathrm{ int}\,G\cap R(L)\ne\varnothing$. It is easy to see that the latter inclusion holds for $L$ and $G$ defined by~\eqref{eq:Lx} and~\eqref{eq:G} respectively. Recalling the definition of $L'$ (formula~\eqref{eq:sLz}) and noting $c^2(G,z-v)=\int_{t_0}^T\langle Q^{-1}(z-v),z-v\rangle dt$ we derive from~\eqref{eq:gerr}-\eqref{eq:supGRL}
$$
\sigma(u)=\min_v\{\int_{t_0}^T(Q^{-1}(z-v),z-v)dt\}+\int_{t_0}^T(R^{-1}u,u)dt
$$ where $L'v=0$. This concludes the proof. 
\end{pf}
\subsection{Optimality conditions and estimation algorithms}
Theorem~\ref{t:1} states that minimax estimation problem is equal to some
 optimal control problem for $\ell\in\mathcal L(T)$, which is called dual control
problem. In the next proposition we
 introduce an approximate solution to the dual control problem without restricting the matrices $F$ and $C(t)$.  
 \begin{prop}\label{p:3e}[Tikhonov regularization]
   Let $\ell\in\mathcal L(T)$. For any $\varepsilon>0$ 
the Euler-Lagrange system
\begin{equation}\label{eq:daebvpe}
\begin{split}
&(F'z)_t(t)=-C'(t)z(t)+H'(t)\hat u+\hat p,\\
&(Fp)_t(t)=C(t)p(t)+\varepsilon Q^{-1}(t)z(t),\\
&\varepsilon\hat u = Rp,Fp(t_0)=0,F'z(T)+F^+Fp(T)=F'\ell
\end{split}
\end{equation} 
has a unique solution $\hat p(\varepsilon),\hat z(\varepsilon)$. 
and 
\begin{equation}
  \label{eq:uze}
  \begin{split}
    &\hat u(\varepsilon):=\frac 1\varepsilon RH\hat p(\varepsilon)\to\hat u\text{ in } 
\mathbb L_2(t_0,T,\mathbb R^p),\\
    &\hat z(\varepsilon)\to\hat z\text{ in }  \mathbb L_2(t_0,T,\mathbb R^m),\\
    &\hat\sigma(\ell)=
\lim_{\varepsilon\to 0}
\frac 1\varepsilon(\langle F'\ell-F^+F\hat p(T),F\hat p(T)\rangle-\int_{t_0}^T\|\hat p(\varepsilon)\|^2dt),
  \end{split}
\end{equation}
where $\hat u,\hat z$ denotes the solution of
\begin{equation}
  \label{eq:dualoc}
  \begin{split}
  &\min_u\{\min_v\{\int_{t_0}^T(Q^{-1}(z-v),z-v)dt\}+\int_{t_0}^T(R^{-1}u,u)dt\},\\
  &(F'z)_t(t)=-C'(t)z(t)+H'(t)u(t), F'z(T)=F'\ell
  \end{split}
\end{equation}
 \end{prop}
 \begin{pf}
   Let $\ell$-error be finite. Then~\eqref{eq:zul} has a solution due to Theorem~\ref{t:1}. Define $(\tilde Hu)=(H'u,0)$, $\tilde l=(0,\ell)$ and set $(Dz)=(-(F'z)_t-C'z,F'z(T))$ for $z(\cdot)\in\mathscr D(D)=\{z(\cdot):F'z(\cdot)\in\mathbb H^1(t_0,T,\mathbb R^n)\}$. It is not difficult to see that the solution to~\eqref{eq:dualoc} coincides with the solution $(\hat u,\hat z)$ of the optimization problem\footnote{The norm is defined by $(u,z)\to\int_{t_0}^T\langle Q^{-1}z,z\rangle+\langle R^{-1}u,u\rangle dt$} $$
\|u\|^2+\|z\|^2\to\min_{u,z}, \quad D'z+\tilde Hu = \tilde l
\eqno(*)$$ This observation allows to apply the Tikhonov regularization~\cite{Tikhonov1977} method in order to derive~\eqref{eq:uze}. For simplicity assume that $Q$ and $R$ are equal to the identity mapping. Let us introduce Tikhonov function
\begin{equation}\label{eq:Tf}
  \begin{split}
    T_\varepsilon&(u,z):=\|F'z(T)-F'\ell\|^2+\int_{t_0}^T \| (F'z)_t+C'z-H'u\|^2dt\\&+\varepsilon\int_{t_0}^T\|u\|^2+\|z\|^2dt 
    =\| Dz-\tilde Hu - \tilde l\|^2+\varepsilon(\|u\|^2+\|z\|^2)
  \end{split}
\end{equation}
It is strictly convex and coercive. 
Thus it's minimum is attained at the unique point $(\hat u(\varepsilon),\hat z(\varepsilon))$. Moreover, $(\hat u(\varepsilon),\hat z(\varepsilon))$ goes to $(\hat u,\hat z)$ in $\mathbb L_2(t_0,T)$ as it follows from properties of the Tikhonov function~\cite{Zhuk2007}. To conclude the proof it is sufficient to show that $(\hat u(\varepsilon),\hat z(\varepsilon))$ verifies~\eqref{eq:daebvpe}. Using the argument of~\cite{Zhuk2007} we derive the Euler-Lagrange equation for $(\hat u(\varepsilon),\hat z(\varepsilon))$:  
\begin{equation}
  \label{eq:eld}
  \begin{split}
    &Dz+\tilde Hu+\tilde p=\tilde l,\\
    &D'\tilde p=\varepsilon z,\quad
    \tilde H'p=\varepsilon u
  \end{split}
\end{equation}
where $\tilde p=(p,q)$, 
$\tilde H'\tilde p=(H'p,0)$, $D'$ is defined by the rule $D'\tilde p=(Fp)_t-Cp$ with $
\tilde p\in\mathscr D(D')=\{\tilde p=(p,q):Fp\in\mathbb H^1(t_0,T,\mathbb R^m),Fp(t_0)=0,q=F^+Fp(T)+d, Fd=0\}
$. For the detailed derivation of $D'$ we refer the reader to~\cite{Zhuk2007}. Now, introducing the definitions of $D$, $D'$, $H$ and $H'$ into~\eqref{eq:eld} we obtain~\eqref{eq:daebvpe}. This proves the existence and uniqueness. We conclude with proving the last line in~\eqref{eq:uze}, which follows from~\eqref{eq:daebvpe} and the formula $$
\|\hat u(\varepsilon)\|^2+\|z(\varepsilon)\|^2\to
\|\hat u\|^2+\|z\|^2=\hat\sigma(\ell)
$$   
 
 \end{pf}
 \begin{rem}
   In fact, the above Proposition claims that $\ell$-estimation $\hat u$ is approximated by $\hat u(\varepsilon)$ for any direction $\ell\in\mathcal L(T)$, provided $\hat u(\varepsilon)$ is a linear transformation of a solution of Euler-Lagrange system~\eqref{eq:daebvpe} for the Tikhonov function~\eqref{eq:Tf} and \eqref{eq:Tf} approximates the minimal worst-case error $\hat\sigma(T,\ell)$. 
 \end{rem}
Now we will derive the suboptimal worst-case recursive estimator, acting\footnote{giving the estimation of the projection of the state vector $x(t)$ onto a minimax observable subspace $\mathcal L(t)$ for all $t\in[t_0,T]$.} on a minimax observable subspace. To do so we will introduce a splitting of~\eqref{eq:daebvpe} into differential and algebraic parts. \\
Let $D=\mathrm{diag}(\lambda_1\dots\lambda_r)$ where $\lambda_i$, $i=\overline{1,r:=\mathrm{rang} F}$ are positive eigen values of $FF'$ and set $\Lambda=\bigl
  (\begin{smallmatrix}
    D^\frac12&&0_{r\times n-r}\\
    0_{m-r\times r}&&0_{m-r\times n-r}
  \end{smallmatrix}\bigr)$. Then \cite{Albert1972} there exist $S_L\in\mathbb R^{m\times m}$, $S_R\in\mathbb R^{n\times n}$ such that 
\begin{equation}
    \label{eq:svd}
      F=S_L\Lambda S_R,S_LS_L'=I, S_RS'_R=I,
 \end{equation} 
Transforming~\eqref{eq:dae} according to~\eqref{eq:svd} and changing the variables one can reduce the general case to the case 
$F=(
\begin{smallmatrix}
  I&&0\\0&&0
\end{smallmatrix}
)$. We split $C(t)$, $Q(t)$ and $H'(t)R(t)H(t)$ according to the structure of $F$ as follows: $
C(t)=\bigl(
\begin{smallmatrix}
C_1&&C_2\\
C_3&&C_4  
\end{smallmatrix}
\bigr)$, $Q=\bigl(
\begin{smallmatrix}
Q_1&&Q_2\\
Q'_2&&Q_4  
\end{smallmatrix}
\bigr)$, $
H'RH=\bigl(
\begin{smallmatrix}
S_1&&S_2\\
S_2'&&S_4
\end{smallmatrix}
\bigr)
$. Define 
\begin{equation*}
  \begin{split}
&A(t)=(C_3'Q_4^{-1}C_4+S_2),\,B(t)=(C'_2-C_4'Q_4^{-1}Q_2'),\\
&C(t,\varepsilon)=-C_1'+C_3'Q_4^{-1}Q_2'+A(t)M(t,\varepsilon)B(t),\\
&W(t,\varepsilon)=(\varepsilon I+S_4+C_4'Q_4^{-1}C_4), M(t,\varepsilon)=
W^{-1}(t,\varepsilon) \\
&Q(t,\varepsilon)=-\frac 1\varepsilon A(t)M(t,\varepsilon)A'(t)+I+\frac 1\varepsilon[S_1+C_3'Q_4^{-1}C_3],\\
&S(t,\varepsilon)=\varepsilon Q_1 -\varepsilon Q_2Q_4^{-1}Q_2'
+\varepsilon B'(t)M(t,\varepsilon)B(t),\\
  \end{split}
\end{equation*}
 \begin{cor}\label{c:eme}[suboptimal estimation on a subspace]
Let
\begin{equation}\label{eq:xKz}
  \begin{split}
   &\hat x_t = -C'(\varepsilon,t)-K(\varepsilon{},t)Q(\varepsilon,t))\hat x + 
   \frac 1\varepsilon\Phi H'R y, \\
   &\dot z_1=C(\varepsilon,t)z_1+Q(\varepsilon,t)Kz_1,z_1(T)=(I+K)^{-1}\ell_1\\
   &K_t=-KC(\varepsilon,t)-C'(\varepsilon,t)K-
   KQ(\varepsilon,t)K\\
   &+S(\varepsilon,t),
   K(0)=0,\hat x(0)=0
    \end{split}
  \end{equation}
with\footnote{$\ell$ is splitted into parts according to the splitting of $\mathbb R^n$ induced by the block structure of $F$.} $\ell=(\ell_1,\ell_2)$ and $\Phi(t,\varepsilon)= \bigl(
\begin{smallmatrix}
  K(t,\varepsilon)\\
M(t,\varepsilon)[\varepsilon B(t)-A'(t)K]
\end{smallmatrix}\bigr)$. Then
\begin{equation*}
  \begin{split}
    \sup_{x_0,f,\eta} &E[\langle\ell,Fx(T)\rangle-\langle (I+K(T,\varepsilon))^{-1}\ell_1,\hat x(T,\varepsilon)\rangle]^2\to \\
    &\inf_u\sup_{x_0,f,\eta} E[\ell(x)-u(y)]^2,\\
    \hat\sigma(\ell)&=\inf_{\varepsilon>0}\frac1\varepsilon\bigr[\langle(I+K(T,\varepsilon))^{-1}\ell_1,K(T,\varepsilon)(I+K(T,\varepsilon))^{-1}
\ell_1\rangle\\
&-\int_{t_0}^T\|\Phi(t,\varepsilon)z_1\|^2dt\bigl]
  \end{split}
\end{equation*}
 \end{cor}
 \begin{pf}
   The idea\footnote{The same idea was used in~\cite{Eremenko1980}.} is to split the Euler-Lagrange system~\eqref{eq:daebvpe} 
into differential $(p_1,z_1)$ and algebraic $(p_2,z_2)$ parts using the splittings of $F$, $C$, $Q$ and $H'RH$, introduced above. We have  
\begin{equation}
  \begin{split}
    &\dot p_1=C_1p_1+C_2p_2+\varepsilon (Q_1z_1+Q_2z_2), p_1(t_0)=0,\\
    &\dot z_1=-C_1'z_1-C_3'z_2+p_1+\frac 1\varepsilon (S_1p_1+S_2p_2),\\
    &0 = C_3p_1+C_4p_2+\varepsilon (Q_2'z_1+Q_4z_2),\\
    &0= -C_2'z_1-C'_4z_2+\frac 1\varepsilon S_2'p_1+(I+\frac 1\varepsilon S_4)p_2\\
    &z_1(T)+p_1(T)=\ell_1
  \end{split}
\end{equation} 
Solving the algebraic equations for $(p_2,z_2)$ 
\begin{equation}\label{eq:p2z2}
    \begin{split}
      &z_2= Q_4^{-1}[(-Q_2'-C_4M)Bz_1+\frac 1\varepsilon (C_4MA'-C_3)p_1],\\
      &p_2=\varepsilon MBz_1- MA'p_1,
    \end{split}
  \end{equation}
and substituting the resulting expressions to into differential equations for $(p_1,z_1)$ one obtains
\begin{equation}
\label{eq:ebvp}
  \begin{split}
    &\dot z_1=C(\varepsilon,t)z_1+Q(\varepsilon,t)p_1, z_1(T)+p_1(T)=\ell_1,\\
    &\dot p_1=-C'(\varepsilon,t)p_1+S(\varepsilon,t)z_1, p_1(t_0)=0
  \end{split}
\end{equation}
Applying simple matrix manipulations one can prove that $Q(\varepsilon,t)\ge0$ and $S(\varepsilon,t)\ge0$ for $\varepsilon>0$, implying \eqref{eq:ebvp} is a non-negative Hamilton system for any $\varepsilon>0$. Therefore it is always solvable and the Riccati equation~\eqref{eq:xKz} has a unique symmetric non-negative solution. Note, that the unique solvability of~\eqref{eq:ebvp} is also implied by Tikhonov method: \eqref{eq:ebvp} is equivalent to the Euler-Lagrange system~\eqref{eq:daebvpe}, which is uniquely solvable. 
Now, by direct calculation  
we derive from~\eqref{eq:p2z2}-\eqref{eq:ebvp}
$$
\hat p=(p_1,p_2)^T = \Phi(t,\varepsilon)z_1
$$
Recalling that (Proposition~\ref{p:3e}) 
$$
\int_{t_0}^T\langle y,\frac 1\varepsilon RH\hat p(\varepsilon)\rangle dt \to
\int_{t_0}^T\langle y,\hat u\rangle\text{ in } \mathbb L_2(t_0,T,\mathbb R^p),
\varepsilon\to 0,
$$ and $\hat p=\Phi(t,\varepsilon) z_1$, $\dot z_1=C(\varepsilon,t)z_1+Q(\varepsilon,t)Kz_1, z_1(T)=(I+K)^{-1}\ell_1$, we derive, integrating by parts, that 
\begin{equation*}
  \begin{split}
    &\int_{t_0}^T\langle y,\frac 1\varepsilon RH\hat p(\varepsilon)\rangle dt = 
\int_{t_0}^T\langle \frac 1\varepsilon \Phi' H'R y, z_1\rangle dt=\\
&\langle (I+K)^{-1}\ell_1,\hat x(T,\varepsilon)\rangle
  \end{split}
\end{equation*}
where $\hat x_t$ is defined in \eqref{eq:xKz}. 
In the same manner we derive the expression for minimax error, recalling \eqref{eq:uze}. 
This concludes the proof. 
 \end{pf}
Now we consider one special case when the DAE is regular and there is a possibility to derive the optimal state estimation algorithm. Let $P^2=P$ ($V^2=V$) and $R(P)=R(F)$ ($R(V)=R(F')$). 
\begin{prop}\label{p:2}[$\ell$-estimation and error]
Let $R((I-V)C'P)\subseteq R((I-V)C'(I-P))$. Then for any $\ell\in\mathbb R^n$  
\begin{equation}
    \label{eq:dae_bvp}
    \begin{split}
      &(Fp)_t(t)=C(t)p(t)+Q^{-1}(t)z(t),Fp(t_0)=0,\\
      &(F'z)_t(t)=-C'(t)z(t)+H'(t)R(t)H(t)p(t),F'z(T)=F'\ell
    \end{split}
  \end{equation}
has a solution. If $p(\cdot)$ and $z(\cdot)$ are some solution of~\eqref{eq:dae_bvp} then, the $\ell$-estimation $\hat u$ is given by $\hat u=RHp$ and the $\ell$-error is represented by $\hat\sigma(T,\ell)=\langle F'^{+}\ell,Fp(T)\rangle$. 
\end{prop}
\begin{pf}
As above~\eqref{eq:svd} we split $C(t)$ and $H'(t)R(t)H(t)$ and $Q$ according to the structure of $F$. 
In this case~\eqref{eq:dae_bvp} reads as 
\begin{equation}
    \label{eq:ade_bvp}
    \begin{split}
      &\dot p_1=C_1p_1+Q_1z_1+C_2p_2+Q_2z_2,p_1(t_0)=0,\\
      &\dot z_1=-C_1'z_1+S_1p_1-C_3'z_2+S_2p_2,z_1(T)=\ell_1,\\
      &0=C_3 p_1+C_4p_2+Q'_2z_1+Q_4z_2,\\
      &0=-C_2'z_1-C_4'z_2+S'_2p_1+S_4p_2
    \end{split}
  \end{equation}
Since
$Q^{-1}>0$ it follows that $Q_4>0$ implying $
z_2=-Q_4^{-1}(C_3p_1+C_4p_2+Q'_2z_1)$ so that
\begin{equation}
  \label{eq:S4p2}
W(t,0)p_2=B(t)z_1-A'(t)p_1
\end{equation} 
where $A,B, W$ were defined above.  
It is easy to see that for our choice of $F$ the proposition's assumption implies $R(C_4^T)\subset R(C_2^T)$. Therefore \eqref{eq:S4p2} is always solvable (in the algebraic sense) and one solution has the form $$
p_2=W^+(t,0)(B(t)z_1-A'(t)p_1)
$$ 
Now we have to assume that 
$p_2\in\mathbb L_2(t_0,T)$. 
Substituting the representation for $p_2$ into \eqref{eq:ade_bvp}
and noting that $C_4(I-W^+(t,0)W(t,0))=0$ we obtain 
\begin{equation}
\label{ade_gmlt}
  \begin{split}
      &\dot p_1=C_+(t)p_1+S_+(t)z_1,z_1(T)=\ell_1 \\
      &\dot z_1=-C_+'(t)z_1+Q_+(t)p_1,p_1(t_0)=0\\
      &
p_1(t_0)=0,z_1(T)=\ell_1
        \end{split}
\end{equation}
where $C_+(t):=C_1-Q_2Q_4^{-1}C_3-B'W^+(t,0)A'$,
$S_+(t):=Q_1-Q_2Q_1^{-1}Q_3+B'W^+(t,0)B$,
$Q_+(t):=S_1+C'_3Q_4^{-1}C_3-AW^+(t,0)A'$.
Applying simple matrix manipulations one can prove that $S_+\ge0$, $Q_+\ge0$ so that \eqref{ade_gmlt} is a non-negative Hamilton system. Therefore it is always solvable. 

With help of~\eqref{eq:dae_bvp} 
one easily shows $I(u)-I(\hat u)\ge 0$ and $I(\hat u)=\hat\sigma=\langle F'^{+}\ell,Fp(T)\rangle$. 

\end{pf}
\begin{rem}
  It is interesting to note that $\frac 1\varepsilon S(t,\varepsilon)\to S_+(t)$, $\varepsilon Q(t,\varepsilon)\to Q_+(t)$ and $-C'(t,\varepsilon)\to C_+(t) $, provided $\varepsilon\downarrow 0$ and the assumptions of the Proposition~\ref{p:2} hold. 
\end{rem}
\begin{cor}\label{c:me}[minimax estimation on a subspace]
  Let $$
\dot K=C_+(t)K+KC_+(t)'-KQ_+(t)K+S_+(t),K(t_0)=0
$$ where $C_+$, $Q_+$ and $S_+$ are defined above (proof of the Proposition~\ref{p:2}). Then $$
\int_{t_0}^T\langle\hat u,y\rangle = 
\langle\ell_1,\hat x(T)\rangle,\hat\sigma(T,\ell)=\langle K(T)\ell_1,\ell_1\rangle
$$ where $\hat x(t_0)=0 $ and 
\begin{equation*}
  \label{eq:flt}
  \hat x_t = (C_+(t)-KQ_+(t))\hat x + 
(\begin{smallmatrix}
  K&&(B'-KA)W^+(t,0)
\end{smallmatrix})
HRy(t)
\end{equation*}
\end{cor}
\begin{pf}
  By direct calculation one finds that $Kz_1,z_1$ verify~\eqref{ade_gmlt} so that $p_1=Kz_1$. Combining this and \eqref{eq:S4p2} with $\hat u=RHp$ and $\hat\sigma(T,\ell)=\langle F'^{+}\ell,Fp(T)\rangle$ (Proposition~\ref{p:2}) one obtains the statement of the corollary. 
\end{pf}
\subsection{Numerical example: non-causal non-stationary DAE}
Let $$
F=\bigl(
\begin{smallmatrix}
  1&&0\\0&&0
\end{smallmatrix}
\bigr), C(t)=\bigl(
\begin{smallmatrix}
  -1&&1\\c_3(t)&&0
\end{smallmatrix}
\bigr) ,H(t)=\bigl(
 \begin{smallmatrix}
   0&&1
 \end{smallmatrix}
 \bigr)
$$
Then $\mathrm{det}(F-\lambda C(t))\equiv 0$ if $c_3(t)=0$. 
The corresponding DAE reads
\begin{equation}
    \label{eq:dae:exmpl}
    \begin{split}
      &\dot x_1 = -x_1 + x_2 + f_1(t),\\
      &0 = c_3(t)x_1(t)+f_2(t), x_1(0)=0
    \end{split}
   \end{equation}
Set $f_1=0$ for simplicity. We have $
x_1(t)=\int_0^t\exp(s-t)x_2(s)ds
$ and 
\begin{equation}
  \label{eq:ineq}
  c_3(t)\int_0^t\exp(s-t)x_2(s)ds = -f_2(t)
\end{equation} 
Set $c_3^+(t)=0$ if $c_3(t)=0$ and $\frac 1{c_3(t)}$ otherwise. Then, formally  
\begin{equation}
  \label{eq:x2}
x_2(t)=
\exp(-t)\dfrac d{dt}(-\exp(t)c^+_3(t)f_2(t))+v(t),
\end{equation}
with $c_3(t)\int_0^t\exp(t-s)v(s)ds=0$ for $ t\in[0,T]$. 
We see that $f_2$ must be able to suppress the growth of $c^+_3$ near points where $c_3(t)=0$, in order to~\eqref{eq:x2} belong to $\mathbb L_2(0,T)$. Taking $f_2=\exp(-c_3^+)b,b(\cdot)\in\mathbb H^1(t_0,T)\Rightarrow $ we obtain $
f_2(t)\in \mathcal R(c_3(t))$ and $c_3^+f_2\in\mathbb H^1(0,T)$. Therefore $x_1(t) = \int_0^t\exp(s-t)x_2(s)ds$ and 
$x_2(t)=-c_3^+(t)\exp(-c_3^+(t))(b(t)+b_t(t))-(c_3^+(t)\exp(-c_3^+(t)))_tb(t)+v(t)$. From this formula we see that $x_2$ is driven by $v$ only if $c_3(t)=0$. If $c_3(t)\ne 0$ then $x_2$ is driven by $f_2$, its derivative and $v$. From the analytical point of view this implies the corresponding DAE is ill-posed: $x_1$ is non-unique and is not continuous with respect to input data. Namely, as the differential operator is unbounded in $\mathbb L_2(t_0,T)$, we see that $x_2$ is not continuous with respect to $f_2$, implying ill-posedness. 
As $x_2$ depends on an arbitrary function $v$ from some linear subspace we have non-uniqueness. 

\begin{figure}
  \centering
\includegraphics[height=.3\textheight]{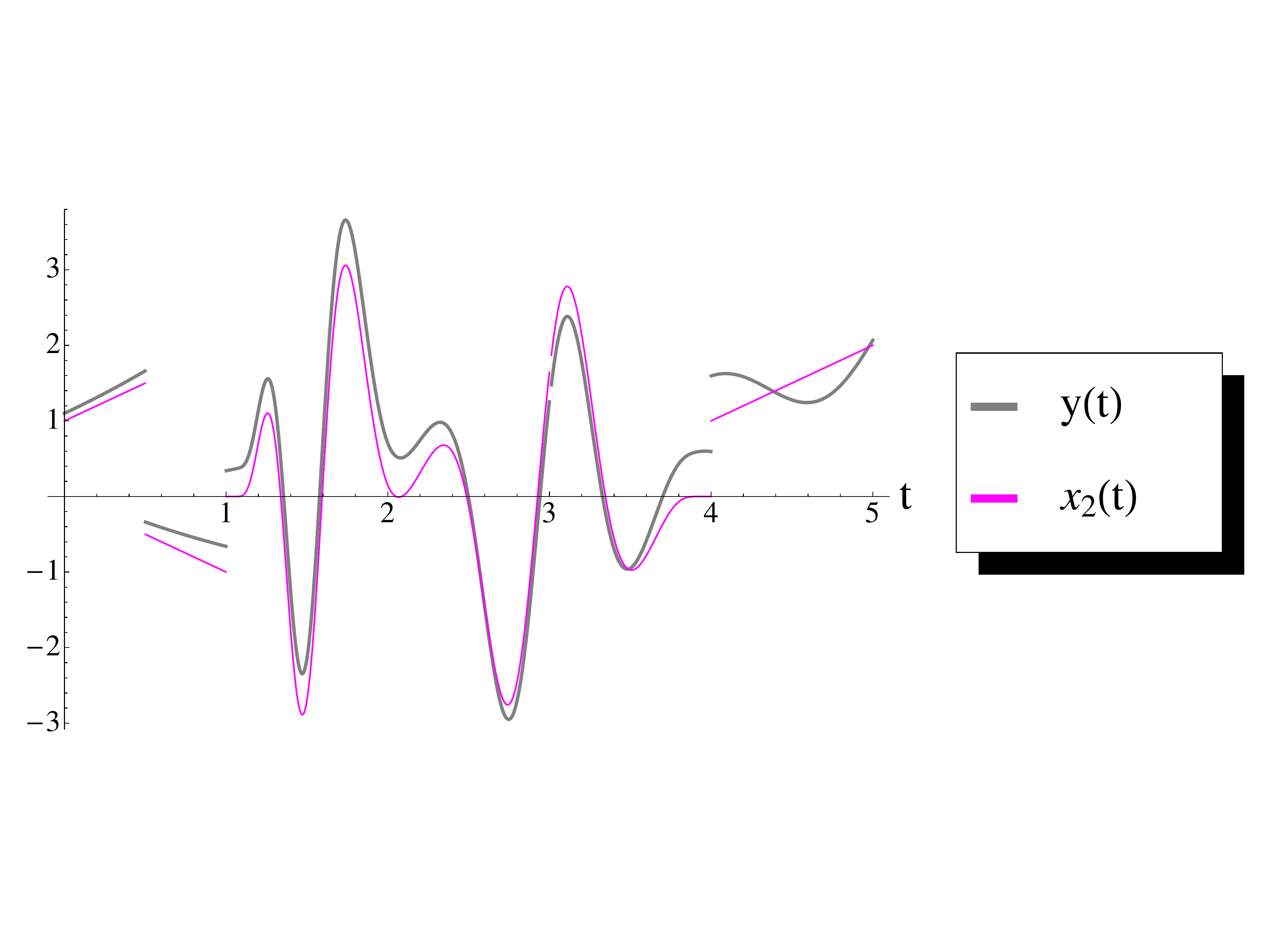}  
  \caption{Observation's noise}
  \label{fig:1}
\end{figure}

\begin{figure}
   \centering
\includegraphics[height=.3\textheight]{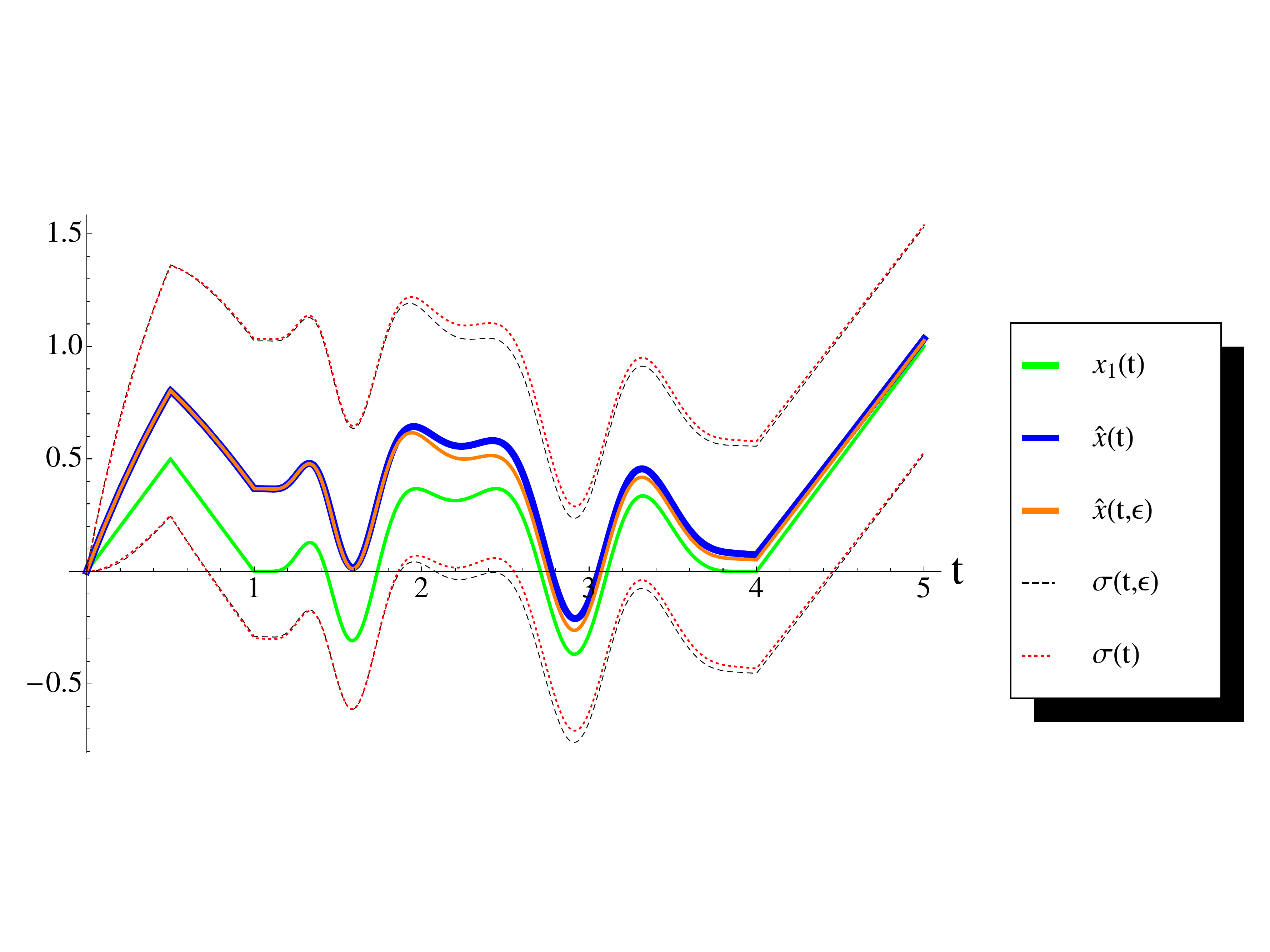}   
   \caption{Optimal estimation and error, suboptimal estimation and error and $x_1$}
   \label{fig:2}
 \end{figure}

The aim is to estimate $x_1(t)$, provided $y(t)=x_2(t)+\eta(t)$ is measured and 
\begin{equation}
  \label{eq:Gn}
\mathrm E\int_{0}^T \frac 6T\eta^2(t)dt\le 1,\quad \mathrm E\eta(\cdot)=0  
\end{equation}
and $(x_1,x_2)$ obeys 
\begin{equation}
    \begin{split}
      &\dot x_1 = -x_1 + x_2+f_1,\quad x_1(0)=0\\
      &0 = c_3(t)x_1(t)+f_2(t)
    \end{split}
   \end{equation}
with 
$f_2(t)=\exp(-c_3^+(t))b(t)$, $b(\cdot)\in\mathbb H^1(0,T)$ and 
\begin{equation}
  \label{eq:Gf}
\int_0^T f_1^2+\frac{\exp(\sqrt{t})}2f_2^2 dt\le 1  
\end{equation}
Figure~\ref{fig:1} describes the observations $y$ of $x_2$, perturbed by the non-Gaussian noise $\eta(\cdot)$, arbitrary function $v(\cdot)$ and uncertain $f_2(\cdot)$, provided $\eta(\cdot)$ verifies~\eqref{eq:Gn}, $f_2(\cdot)$ verifies~\eqref{eq:Gf} and $c_3(t)\int_0^t\exp(t-s)v(s)ds=0$ for $ t\in[0,T]$. 
The sub-optimal estimation $\hat x(t,\varepsilon)$ and error are given by
      \begin{equation*}
     \begin{split}
      &\hat x_t(t,\varepsilon) = \bigr[-1-K(1+\frac{2c_3^2(t)}{\varepsilon \exp(\sqrt{t})})\bigl]\hat x(t,\varepsilon) + 
 \frac {6 y(t)}{6+T\varepsilon},\hat x(\varepsilon,0)=0,\\
    &K_t=-2K-(1+\frac{2c_3^2(t)}{\varepsilon \exp(\sqrt{t})})K^2+
  \varepsilon(1+\frac 1{\frac 6T+\varepsilon}),
  \, K(0,\varepsilon)=0,\\
&\hat\sigma(t,\varepsilon):=
\frac 1\varepsilon\bigr[\frac{K(t,\varepsilon)}{(K(t,\varepsilon)+1)^2}-
\int_{0}^t[K^2(s,\varepsilon)+\frac{\varepsilon^2}{(\frac 6T+\varepsilon)^2}]z^{\varepsilon}(s;t)^2ds\bigl],\\
&z_t^{\varepsilon}(t;\tau)=[1+K(1+\frac{2c_3^2(t)}{\varepsilon \exp(\sqrt{t})})]z^{\varepsilon}(t;\tau),
z^{\varepsilon}(\tau,\tau)=(K(\tau,\varepsilon)+1)^{-1}
\end{split}
 \end{equation*}
Note that the pencil $sF'-C'-H'H$ is regular. As~\eqref{eq:dae_bvp} is solvable in this case we apply Proposition~\ref{p:2} in order to derive the minimax estimation $\hat x$ and minimax estimation error $\hat\sigma(t)$:
      \begin{equation*}
     \begin{split}
      &\hat x_t = \bigr[-1-K(t)\frac{2c_3^2(t)}{\exp(\sqrt{t})})\bigl]\hat x + 
  y(t),\quad\hat x(0)=0,\\
    &K_t=-2K-\frac{2c_3^2(t)}{\exp(\sqrt{t})})K^2+
  (1+\frac T{ 6}),
  \, K(0,\varepsilon)=0,\\
&\hat\sigma^2(t,\ell_1):=
K(t)\ell_1^2
     \end{split}
 \end{equation*}
In Figure~\ref{fig:2} the comparison of the optimal estimator $\hat x$ and error $\hat\sigma(t)$ with sub-optimal estimator $\hat x(\varepsilon)$ and error $\hat\sigma(t,\varepsilon)$ are presented, provided $\varepsilon=\exp(-100)$ and $\ell_1=0$. 
\bibliography{refs,myrefs}
\end{document}